%%a Plain TeX file by Shalosh B. Ekhad (2 pages)

%begin macros

\def\N{{\cal N}}
\def\Z{{\cal Z}}
\baselineskip=14pt
\parskip=10pt

\font\eightrm=cmr8 
\font\eighttt=cmtt8
\magnification=\magstephalf

\def\1{{\overline{1}}}
\def\2{{\overline{2}}}
\parindent=0pt
\overfullrule=0in
\def\Tilde{\char126\relax}
\def\frac#1#2{{#1 \over #2}}
%\headline={\rm  \ifodd\pageno  \RightHead  \else  \LeftHead  \fi}
%\def\RightHead{\centerline{
%Title
%}}
%\def\LeftHead{ \centerline{Doron Zeilberger}}
%end macros
\bf
\centerline
{
A Note on the Stanley Distribution
}
\rm
\bigskip
\centerline{ {\it
Shalosh B.
EKHAD}\footnote{$^1$}
{\eightrm  \raggedright
Department of Mathematics, Rutgers University (New Brunswick),
Hill Center-Busch Campus, 110 Frelinghuysen Rd., Piscataway,
NJ 08854-8019, USA.
%\break
c/o {\eighttt zeilberg  at math dot rutgers dot edu} ,
\hfill \break
{\eighttt http://www.math.rutgers.edu/\~{}zeilberg/} .
Jan. 20, 2010. Exclusively published in the Personal Journal of Ekhad and Zeilberger,
{\eighttt http://www.math.rutgers.edu/\~{}zeilberg/pj.html} , and {\tt arxiv.org} \quad .
}
}
In the second of the three AMS colloquium talks, recently masterfully\footnote{$^2$}
{\eightrm  \raggedright
The three colloquium talks delivered each year during the annual meeting of the American Mathematical Society
are supposed to be the highlights of the meeting, and, in the old days, every participant was expected to attend,
since they were supposed (often successfully) to explain to non-specialists, in a lucid and engaging manner, some non-trivial
mathematics, and were delivered by a great mathematician. Some of these series of talks were better than others, and the three lectures
by Richard Stanley, last week, were amongst the very best ones
(along with Gian-Carlo Rota's 12 years ago, and a few other ones). Unfortunately, only about
400 participants (out of more than 5000 people!) attended the first lecture (where there was no conflict with
shorter talks), and only about 150 attended the second and third ones. Shame on you,
2010 ``mathematicians''!
\hfill \break
Luckily, one can still retrieve the slides of these fascinating talks from Richard Stanley's website:
{\eighttt http://www-math.mit.edu/\Tilde rstan/transparencies/ids-ams.pdf} ,
{ \eighttt  http://www-math.mit.edu/\Tilde rstan/transparencies/altperm-ams.pdf} ,
{\eighttt http://www-math.mit.edu/\Tilde rstan/transparencies/redec-ams.pdf } .
}
presented by Algebraic Combinatorics guru Richard Stanley, he talked about the
random variable ``length of the largest alternating (i.e. up-down) subsequence'' in a random
permutation of length $n$. This was in analogy with the celebrated random variable
``length of the largest increasing subsequence''. For the latter, it was famously proved by
Baik-Deift-Johansson that the {\it limiting distribution} (after it is centralized and
divided by the standard deviation) is the intriguing Tracy-Widom distribution.

Stanley humorously narrated that when he investigated the limiting distribution for the former random variable,
he was hoping that the limiting distribution would be equally interesting, and has already fantasized
that he would be immortalized (along with Gauss, Poisson, Cauchy, and Tracy-Widom) by having
an exotic new probability distribution named after him. To his dismay, it turned out
(because of a general theorem of Robin Pemantle and Herb Wilf) that the limiting distribution
is the {\it utterly boring} {\bf Gaussian} (aka {\it normal}) distribution. So much for Stanley's
dream of immortality (of course, there is a {\it Stanley-Reisner Ring}, but ``{\it Stanley Distribution}''
has a better {\it ring} to it!).

Stanley found that the
{\it expectation} $\mu_n$ and the {\it variance} $\sigma_n^2$ are given by
$$
\mu_n= \frac{2}{3}\,n \, + \, \frac{1}{6} \quad ; \quad
\sigma_n^2 ={\frac {8}{45}}\,n \, -{\frac {13}{180}}\quad  (for \quad n  \geq 4 \,) \quad ,
$$
and he deduced (from the explicit generating function that he derived, using the above-mentioned Pemantle-Wilf theorem) that
$$
\Z_n := \frac{X_n-\mu_n}{\sigma_n} 
$$
converge to the normal distribution, i.e.:
$$
\Z_n \rightarrow \N \quad (in \quad distribution) .
$$
This is equivalent to the statement that the moments of $\Z_n$ converge to the moments of the Gaussian distribution
($0$ for odd moments and $1 \cdot 3 \cdot \dots \cdot (2r-1)= \frac{(2r)!}{2^r r!}$ for the even, $2r$-th, moment),
to wit:
$$
\alpha_{2r} (\Z_n) = \frac{(2r)!}{2^r r!}+ o(1) \quad ,
$$
$$
\alpha_{2r+1} (\Z_n) =  o(1) \quad .
$$
So the leading terms for the asymptotics for the moments are indeed the {\it boring}, normal ones. {\bf But},
thanks to the amazing Maple package {\tt HISTABRUT}, written by my master, Doron Zeilberger,
that will soon be released with an accompanying article explaining how to use it,
one can get:
$$
\alpha_{2r}(\Z_n)=
{{(2r)!} \over {2^r r!}}
\, 
\left( 1+{\frac {1}{1764}}\,{\frac {r \left( r-1 \right)  \left( 10\,r-713 \right) }{n}} + O({{1} \over {n^2}}) \right)  
 \quad ,
$$
$$
\alpha_{2r+1}(\Z_n)=
$$
$$
-{\frac {\sqrt{10}} {43}}\, 
{{(2r)!} \over {2^r r!}}
{\frac {1}{\sqrt {n}}} 
\left( r-1+{\frac {1}{931392}}\,{\frac { \left( r-1 \right)  \left( 1760\,{r}^{3}-381744\,{r
}^{2}+1430752\,r+150351 \right) }{n}} \, + \, O({{1} \over {n^2}}) \right)   \quad .
$$
These formulas tell you {\it how fast} the moments of the (discrete) Stanley Distribution ($\Z_n$) converge to 
those of the  normal distribution.

Zeilberger's Maple package {\tt HISTABRUT} can yield even higher-order asymptotics, but {\it who cares?}.
Also, I have to admit that the above formulas are not yet rigorous, but they are certainly {\it rigorizable}.
So there exists a rigorous proof, but {\it who cares}?

\end